\documentclass{amsart}
\usepackage{amsmath}
\usepackage{graphicx}
\usepackage{amsfonts}
\usepackage{amssymb}
\newtheorem{theorem}{Theorem}
\theoremstyle{plain}

\newtheorem{definition}{Definition}

\newtheorem{proposition}{Proposition}
\newtheorem{remark}{Remark}

\numberwithin{equation}{section}
\begin{document}
\title{On the Betti numbers of a loop space}
\author{Samson Saneblidze}
\address{A. Razmadze Mathematical Institute\\
Department of Geometry and Topology\\
M. Aleksidze st., 1\\
0193 Tbilisi, Georgia} \email{sane@rmi.acnet.ge} \dedicatory{To Tornike Kadeishvili and
Mamuka Jibladze}
\date{}

\begin{abstract}
Let $A$ be a special homotopy G-algebra over a commutative unital ring $\Bbbk$ such
that both $H(A)$ and $\operatorname{Tor}_{i}^{A}(\Bbbk,\Bbbk)$ are finitely generated
$\Bbbk$-modules for all $i$, and let $\tau_{i}(A)$ be the cardinality of a minimal
generating set for the $\Bbbk$-module $\operatorname{Tor}_{i}^{A}(\Bbbk,\Bbbk).$ Then
the set $\left\{  \tau _{i}(A)\right\}  $ is unbounded if and only if $\tilde{H}(A)$
has two or more algebra generators. When $A=C^{\ast}(X;\Bbbk)$ is the simplicial
cochain complex of a simply connected finite $CW$-complex $X,$ there is a similar
statement for the "Betti numbers" of the loop space $\Omega X.$ This unifies existing
proofs over a field $\Bbbk$ of zero or positive characteristic.

\end{abstract}
\maketitle

\section{Introduction}

Let $Y$ be a topological space, let $\Bbbk$ be a commutative ring with identity, and
assume that the $i^{th}$-cohomology group $H^{i}(Y;\Bbbk)$ of $Y$ is finitely generated
as a $\Bbbk$-module. We refer to the cardinality of a minimal generating set of
$H^{i}(Y;\Bbbk),$ denoted by $\beta_{i}(Y),$ as the \emph{generalized
}$i^{th}$\emph{-Betti number} \emph{of} $Y.$

\begin{theorem}
\label{betti} Let $X$ be a simply connected space. If $H^{\ast}(X;\Bbbk)$ is finitely
generated as a $\Bbbk$-module and $H^{\ast}(\Omega X;\Bbbk)$ has finite type, then the
set of generalized $i^{th}$-Betti numbers $\left\{ \beta_{i}(\Omega X;\Bbbk)\right\}  $
is unbounded if and only if $\tilde {H}^{\ast}(X;\Bbbk)$ has at least two algebra
generators.
\end{theorem}

Theorem \ref{betti} was proved by Sullivan \cite{Sul} over fields of characteristic
zero and by McCleary \cite{mcCleary} over\ fields of positive characteristic. However,
Theorem \ref{betti} is a consequence of the following more general algebraic fact: Let
$A^{\prime}=\{{A^{\prime}}^{i}\},i\geq0,$ with ${A^{\prime}}^{0}=\mathbb{Z},$
${A^{\prime}}^{1}=0,$ be a torsion free graded abelian group endowed with a homotopy
$G$-algebra (hga) structure. Then for $A=A^{\prime}\otimes_{\mathbb{Z}}\Bbbk$ we have
the following theorem whose proof appears in Section \ref{theorem}:

\begin{theorem}
\label{tau}Assume that $H^{\ast}(A)$ is finitely generated as a $\Bbbk$-module and that
$Tor_{\ast}^{A}(\Bbbk\,,\Bbbk)$ has finite type. Let $\tau_{i}(A)$
denote the cardinality of a minimal generating set of $Tor_{i}^{A}%
(\Bbbk\,,\Bbbk).$ Then the set $\left\{  \tau_{i}(A)\right\}  $ is unbounded if and
only if $\tilde{H}(A)$ has at least two algebra generators.
\end{theorem}

Let $C^{\ast}(X;\Bbbk)=C^{\ast}(\operatorname{Sing}^{1}X;\Bbbk)/C^{>0}%
(\operatorname{Sing}\,x\,;\Bbbk)$ in which ${\operatorname{Sing}}^{1}%
X\subset{\operatorname{Sing}}X$ is the Eilenberg 1-subcomplex generated by the singular
simplices that send the 1-skeleton of the standard $n$-simplex $\Delta^{n}$ to the base
point $x$ of $X.$ To deduce Theorem \ref{tau} from Theorem \ref{betti}, set
$A=C^{\ast}(X;\Bbbk),$ and apply Proposition
\ref{barVX} below together with the filtered hga model $(RH(A),d_{h}%
)\rightarrow A$ of $A$ (a special case of the filtered Hirsch algebra \cite{sane2}).
Let $BA$ denote the bar construction of $A.$ When $\tilde {H}(A)$ has at least two
algebra generators, we construct two infinite sequences in the filtered model and take
all possible $\smile_{1}$-products of their components to detect a submodule of
$H^{\ast}(BA)$  at least as large as  the polynomial algebra $\Bbbk\lbrack x,y].$

Each of the sequences mentioned above can be thought of as generalizations of an
infinite sequence ($\infty$-\emph{implications} of its first component) introduced by
Browder \cite{browder}. Indeed, this work arose after writing down these special
sequences in the hga resolution of a commutative graded
algebra (cga) over the integers via formulas (\ref{xrelation1}%
)--(\ref{xrelation3}) below, at which point we realized that their construction mimics
that of Massey symmetric products defined by Kraines \cite{kraines} (see also
\cite{sane2}). In general, a sequence formed from Massey symmetric products is closely
related to the one obtained from $A_{\infty}$-operations in an $A_{\infty}$-algebra
defined by Stasheff \cite{Stasheff} by restricting to the same variables in question.
When a differential graded algebra (dga) $A$ is free as a $\Bbbk$-module, the sequence
of $A_{\infty}$-operations on the homology $H\left(  A\right)  $ was constructed by
Kadeishvili \cite{kade}.

I am grateful to Jim Stasheff for comments and suggestions and to the referee for
comments that helped to improve the exposition.

\section{Some preliminaries and conventions}

We adopt the notations and terminology of \cite{sane2}. We fix a ground ring $\Bbbk$
with identity, a primary example of which is the integers $\mathbb{Z}.$ Let
$\mathbb{Z}_{\Bbbk}\subset\mathbb{Z}$ be the subset defined by
\[
\mathbb{Z}_{\Bbbk}=\{\lambda\in\mathbb{Z}\,|\,\lambda_{\Bbbk}:\Bbbk
{\rightarrow}\Bbbk,\,\kappa\rightarrow\lambda\kappa,\ \ \text{is injective}\}.
\]
Let $\mu\in\mathbb{Z}\setminus\mathbb{Z}_{\Bbbk}$ denote the smallest integer such that
$\mu\kappa=0$ for all $\kappa\in\Bbbk.$ Thus if $\mu=0$,
$\mathbb{Z}_{\Bbbk}=\mathbb{Z}\setminus0$ (e.g. $\Bbbk$ is a field of characteristic
zero).

A (positively) graded algebra $A$ is 1-reduced if $A^{0}=\Bbbk$ and $A^{1}=0.$ For a
general definition of an homotopy Gerstenhaber algebra (hga) $\left(
A,d,\cdot\,,\{E_{p,q}\}\right)  _{p\geq0,\,q=0,1}$ see \cite{Gerst-Voron},
\cite{Getz-Jones}, \cite{KScubi}. The defining identities for an hga are the following:
Given $k\geq1,$
\begin{equation}%
\begin{array}
[c]{llll}%
dE_{k,1}(a_{1},...,a_{k};b) & = & \sum_{i=1}^{k}\ \,\,(-1)^{\epsilon
_{i-\!1}^{a}}\,E_{k,1}(a_{1},...,da_{i},...,a_{k};b) & \\
&  & \ \ \ \ +\ \ \,(-1)^{\epsilon_{k}^{a}}\ \ \ E_{k,1}(a_{1},...,a_{k};db) &
\\
&  & \!\!\!\!\!\!+\sum_{i=1}^{k-1}\,\,\,(-1)^{\epsilon_{i}^{a}}\ \ \ E_{k-1,1}%
(a_{1},...,a_{i}a_{i+1},...,a_{k};b) & \\
&  & \ \ \ \ +\ \ \,(-1)^{\epsilon_{k}^{a}+|a_{k}|\!|b|}\,E_{k-1,1}%
(a_{1},...,a_{k-1};b)\!\cdot\!a_{k} & \\
&  & \ \ \ \ +\ \ \,(-1)^{|a_{1}|}\ \ \ \ \ \ \ a_{1}\!\cdot\!E_{k-1,1}%
(a_{2},...,a_{k};b), &
\end{array}
\label{dif}%
\end{equation}%
\begin{multline}
E_{k,1}(a_{1},\!...,a_{k};b\cdot c)\label{prod}\\
=\sum_{i=0}^{k}(-1)^{|b|(\epsilon_{i}^{a}+\epsilon_{k}^{a})}E_{i,1}%
(a_{1},\!...,a_{i};b)\cdot E_{k-i,1}(a_{i+1},\!...,a_{k};c)
\end{multline}
and
\begin{multline}
\sum_{\substack{_{k_{1}+\cdots+k_{p}=k} \\_{1\leq p\leq k+\ell}}%
}(-1)^{\epsilon}E_{p,1}\left(  E_{k_{1},\ell_{1}}(a_{1},...,a_{k_{1}}%
;b_{1}^{\prime}),\!...,E_{k_{p},\ell_{p}}(a_{_{k-k_{p}+1}},...,a_{k}%
;b_{p}^{\prime})\,;c\right) \label{associativity}\\
=E_{k,1}\left(  a_{1},...,a_{k};E_{{\ell},1}(b_{1},...,b_{\ell};c)\right)  ,\\
\hspace{1.4in}b_{i}^{\prime}\in\{1,b_{1},..,b_{\ell}\},\ \ \ \epsilon
=\sum_{i=1}^{p}(|b_{i}^{\prime}|+1)(\varepsilon_{k_{i}}^{a}+\varepsilon
_{k}^{a}),\,b_{i}^{\prime}\neq1,\\
\varepsilon_{i}^{a}=|a_{1}|+\cdots+|a_{i}|+i.
\end{multline}

A \emph{morphism} $f:A\rightarrow A^{\prime}$ of hga's is a dga map $f$ commuting with
all $E_{k,1}.$

\begin{remark}
Note that we do not use axiom (\ref{associativity}) in the sequel.
\end{remark}

Below we review the notion of an hga resolution of a cga as a special Hirsch algebra
(the existence of such a resolution is proved in \cite{sane2}). Given a cga $H,$ its
hga resolution is a multiplicative resolution
\[
\rho:(R^{\ast}H^{\ast},d)\rightarrow H^{\ast},\ \ \ RH=T(V),\ \ \ V=\langle
\mathcal{V}\rangle,
\]
endowed with an hga structure
\[
E_{k,1}:RH^{\otimes k}\otimes RH\rightarrow  RH,\ \ k\geq1,
\]
together with a decomposition of $V$ such that $V^{\ast,\ast}=\mathcal{E}
^{\ast,\ast}\oplus{U}^{\ast,\ast},$ where $\mathcal{E}^{\ast,\ast
}=\{\mathcal{E}_{p,q}^{<0,\ast}\}$ is distinguished by an isomorphism of modules
\[
E_{k,1}:\otimes_{r=1}^{k}R^{i_{r}}H^{k_{r}}\bigotimes V^{j,\ell
}\overset{\approx}{\longrightarrow}\mathcal{E}_{k,1}^{s-k\,,\,t}\subset V^{k-s,t},\ \ \
(s,t)\!=\!\left(  \sum_{r=1}^{k}i_{r}\!+\!j,\sum_{r=1} ^{k}k_{r}\!+\!\ell\right)  .
\]
Furthermore,  if $H$ is a $\mathbb{Z}$-algebra, its hga resolution $(RH,d)$ is
automatically endowed with two operations $\cup_{2}$ and $\smile_{2}$. The first
operation $\cup_{2}$ appears because each cocycle $a\smile_{1}a\in\mathcal{E}_{1,1}\cap
R^{-1}H^{2j},$ where $a\in R^{0}H^{2j},$ is killed by some element in $R^{-2}H^{2j},$
denoted by $a\cup_{2}a.$ The second operation arises from the non-commutativity of
$\smile_{1}$-product in the usual way, and satisfies Steenrod's formula for the
$\smile_{2}$-cochain operation. These two operations are related to each other by the
initial relations $a\smile_{2}a=2a\cup_{2}a$ and $a\smile_{2}b=a\cup_{2}b,$ $a\neq
b\in\mathcal{U}$ with $\langle \mathcal{U}\rangle=U.$  Note also that
$a\smile_{2}a=a\cup_{2}a=0$ for $a\in U$ of odd degree. In general,
$U=\mathcal{T}\oplus \mathcal{N},$ with an element of $\mathcal{T}$ given by
$a_{1}\cup_{2}\cdots\cup_{2}a_{n},$ $a_{i}\in U,\,n\geq2.$ The action of the resolution
differential $d$ on elements of $\mathcal{T}$ such that $da_{i}=0$ is
\begin{multline}
d(a_{1}\cup_{2}\cdots\cup_{2}a_{n})\label{cup2q}\\
=\sum_{(\mathbf{i};\mathbf{j})}(-1)^{|a_{i{_{1}}}\!|+\cdots+|a_{i\!{_{k}}}%
\!|}(a_{i_{1}}\cup_{2}\cdots\cup_{2}\,a_{i_{k}})\,\smile_{1}\,(a_{j_{1}}%
\cup_{2}\cdots\cup_{2}\,a_{j_{\ell}}),
\end{multline}
where we sum over all unshuffles $(\mathbf{i};\mathbf{j})=(i_{1}<\cdots
<i_{k}\,;j_{1}<\cdots<j_{\ell})$ of $\underline{n}$ with $(a_{i_{1}%
},...,a_{i_{k}})=(a_{i_{1}^{\prime}},...,a_{i_{k}^{\prime}})$ if and only if
$\mathbf{i}=\mathbf{i}^{\prime}$ and $\smile_{1}$ denotes $E_{1,1}.$ In particular, for
$a_{1}=\cdots=a_{n}=a=a^{\cup_{2}1}$ and $n\geq2$ we get
$da^{\cup_{2}n}=\sum_{k+\ell=n}a^{\cup_{2}k}\smile_{1}a^{\cup_{2}\ell },\,k,\ell\geq1.$
And in general $d(a\smile_{2}b)=nd(a\cup_{2}b),n\geq1.$

An hga resolution $(RH,d)$ is \emph{minimal} if
\[
d(U)\subset\mathcal{E}+\mathcal{D}+\kappa\!\cdot\!V
\]
where ${\mathcal{D}^{\ast,\ast}}\subset R^{\ast}H^{\ast}$ denotes the submodule of
decomposables $RH^{+}\!\cdot RH^{+}$ and $\kappa\in\Bbbk$ is non-invertible; For
example, $\kappa\in {\mathbb Z}\setminus \{-1,1\}$ when $\Bbbk=\mathbb{Z}$ and
$\kappa=0$ when $\Bbbk$ is a field.

Let   $K=\{K^j\}_{j\geq 3}$ with $K^j=\{a\in \mathcal{V}^{-1,j} \,|\,da=\lambda b,\, \lambda \neq \pm
1,\,b\in \mathcal{V}^{0,j}  \}.$  Note that a general form of a relation in (minimal)
$(RH,d)$ starting by variables $v_i\in {K}\cup\mathcal{V}^{0,*}$ is
\begin{multline}
\label{general} du=\sum_{s\geq1}\lambda_{s} P_{s}(v_{1},...,v_{r_{s}})+\lambda v,\ \ \
\lambda\neq \pm 1,\,\lambda_{s}\neq 0,
\,r_{s}\geq1,  \\
 u\in\bigcup_{i\geq 1}\mathcal{V}^{-i,*},\,\, v\in \bigcup_{i\geq 1}\mathcal{V}^{-i,*}\setminus {K},
\end{multline}
where $P_{s}(v_{1},...,v_{r_{s}})$ is a monomial in $\mathcal{D}^{*,*}\subset R^*H^*.$

 Let $A$ be an hga
and let $\rho:(RH,d)\rightarrow H$ be an hga resolution. A \emph{filtered hga
model} of $A$ is an hga quasi-isomorphism
\[
f:(RH,d_{h})\rightarrow(A,d_{A})
\]
in which
\[
d_{h}=d+h,\ \ \ h=h^{2}+\cdots+h^{r}+\cdots,\ \ \ h^{r}:R^{p}H^{q}\rightarrow
R^{p+r}H^{q-r+1}.
\]
The equality $d_{h}^{2}=0$ implies the sequence of equalities
\[
dh^{2}+h^{2}d=0,\,\,\,\,dh^{3}+h^{3}d=-h^{2}h^{2},\,\,\,\,dh^{4}+h^{4}%
d=-h^{2}h^{3}-h^{3}h^{2},\dotsc,
\]
and $h$ is referred to as a \emph{perturbation of} $d.$ The map $h^{r}%
|_{R^{-r}H}:R^{-r}H\rightarrow R^{0}H,\,r\geq2,$ denoted by $h^{tr},$ is referred to as
the \emph{transgressive} component of $h.$ The fact that the perturbation $h$ acts as a
derivation on elements of $\mathcal{E}$ implies $h^{tr}|_{\mathcal{E}}=0.$ For the
existence of the filtered model see \cite{sane2}.

In the sequel, $A^{\prime}$ denotes a 1-reduced torsion free hga over
$\mathbb{Z}$, while $A$ denotes the tensor product hga $A^{\prime}
\otimes_{\mathbb{Z}}\Bbbk.$ Denote also $H=H^{\ast}(A^{\prime})$ and
$H_{\Bbbk}=H^{\ast}(A).$ Assume $(RH,d)$ is minimal and let $RH_{\Bbbk
}=RH\otimes_{\mathbb{Z}}\Bbbk;$ in particular, $RH_{\Bbbk}=T(V_{\Bbbk})$ for
$V_{\Bbbk}=V\otimes_{\mathbb{Z}}\Bbbk.$ When $\Bbbk$ is a field of characteristic zero,
$\rho\otimes1:RH_{\Bbbk}\rightarrow H\otimes_{\mathbb{Z}}\Bbbk=H_{\Bbbk}$ is an hga resolution of
$H_{\Bbbk},$ which is \emph{not }minimal when $\operatorname{Tor}H\neq0$. In general,
given a filtered model $(RH,d_{h})$ of $A^{\prime},$ we obtain an hga model
\[
f\otimes1:(RH_{\Bbbk},d_{h}\otimes1)\rightarrow(A,d_{A}).
\]
for $(A,d_{A}).$ Denote $\bar{V}_{\Bbbk}=s^{-1}(V_{\Bbbk}^{>0})\oplus\Bbbk$ and define
the differential $\bar{d}_{h}$ on $\bar{V}_{\Bbbk}$ by the restriction of $d+h$ to
$V_{\Bbbk}$ and obtain the cochain complex $(\bar{V}_{\Bbbk },\bar{d}_{h}).$

Since the map $f\otimes1$ is in particular a homology isomorphism (by the universal
coefficient theorem), the following two propositions follow immediately from the
results in \cite{F-H-T} and the standard isomorphisms $H^{\ast}(BA,d_{_{BA}})\approx
Tor^{A}(\Bbbk,\Bbbk)$ and $H^{\ast}(BC^{\ast }(X;\Bbbk),d_{_{BC}})\approx
H^{\ast}(\Omega X;\Bbbk).$

\begin{proposition}
\label{barV} There are isomorphisms
\[
H^{*}(\bar{V}_{\Bbbk},\bar{d}_{h})\approx H^{*}(B(RH_{\Bbbk}),d_{_{B(RH_{\Bbbk })}})
\approx H^{*}(BA,d_{_{BA}}) \approx Tor^{A}(\Bbbk,\Bbbk).
\]

\end{proposition}

And for $A=C^{\ast}(X;\Bbbk)$ we obtain:

\begin{proposition}
\label{barVX} There are isomorphisms
\[
H^{*}(\bar{V}_{\Bbbk},\bar{d}_{h})\approx H^{*}(BC^{*}(X;\Bbbk),d_{_{BC}}) \approx
H^{*}(\Omega X;\Bbbk) .
\]

\end{proposition}

Given $(RH,d)$ and $x,c\in V$ with $dx,dc\in{\mathcal{D}}+\lambda V,\,\lambda\neq1,$
let $\eta_{x,c}$ denote an element of $\mathcal{E}_{>1,1}$ such that
\[
x\smile_{\mathbf{1}}c:=\eta_{x,c}+x\smile_{1}c
\]
satisfies $d(x\smile_{\mathbf{1}}c)\in{\mathcal{D}}+\lambda V.$ For example,
if $dx\in\lambda V,$ then $\eta_{x,c}=0,$ and if $dx=\sum_{i}a_{i}%
b_{i}+\lambda v$ with $da_{i},db_{i}\in\lambda V,$ then $\eta_{x,c}=\sum
_{i}(-1)^{|a_i|}E_{2,1}(a_{i},b_{i}\,;c).$ In general, $\eta_{x,c}$ can be found as
follows: Let $j:B(RH)\rightarrow\overline{RH}\rightarrow\bar{V} $ be the canonical
projection used by the proof of the first isomorphism in Proposition \ref{barV}, and
choose $y\in B(RH)$ so that $j(y)=\bar{x}$ and $j\mu
_{E}(y;\bar{c})=\bar{\eta}_{x,c}+\overline{x\smile_{1}c},$ where the product
$\mu_{E}:B(RH)\otimes B(RH)\rightarrow B(RH)$ is determined by the hga structure on
$RH.$

The following proposition is simple but useful. Let ${\mathcal D}_{\Bbbk}\subset RH$ be
a subset defined by ${\mathcal D}_{\Bbbk}={\mathcal D}$ for $\mu=0$ and
\[
\mathcal{D}_{\Bbbk}=\left\{  u+\lambda v|\,u\in\mathcal{D},v\in V,\,\lambda \ \text{is
divisible by}\ \mu\right\} \ \ \ \text{for}\ \  \mu\geq 2.
\]

\begin{definition}
An element $x\in V$ with $d_{h}x\in{\mathcal{D}}+ \lambda V,\,\lambda\neq1,$ is
\underline{$\lambda$-homologous} \underline{to zero}, denoted by ${[\bar
{x}]}_{\lambda}=0,$ if there are $u,v\in V$ and $z\in\mathcal{D}$ such that
\[
d_{h}u=x+z+\lambda v; \] $x$ is \underline{weakly homologous} to zero when $v=0$ above.
\end{definition}

\begin{proposition}
\label{basic} Let $c\in V$ and $d_{h}c\in\mathcal{D}_{\Bbbk}.$ If $d_{h}c$ has a
summand component $ab\in\mathcal{D}$ such that $a,b\in{V},$ $d_{h}
a,d_{h}b\in\mathcal{D}_{\Bbbk},$ both ${a}$ and  ${b}$  are not weakly homologous to
zero, then  $c$ is also not weakly homologous to zero.
\end{proposition}
\begin{proof}
The proof is straightforward using the equality $d^{2}_{h}=0.$
\end{proof}
In particular,  for $\Bbbk={\mathbb Z},$ under hypotheses of the proposition if
$[\bar{a}],[\bar{b}]\neq0,$ then $[\bar {c}]\neq0$ in $H^{\ast}(\bar{V},\bar{d}_{h}).$

Note that over a field $\Bbbk$, Proposition \ref{basic} reflects the obvious fact that
$x\!\in\! H^{\ast}(\Omega X;\Bbbk)$ is non-zero whenever
 some $x^{\prime}\otimes x^{\prime\prime}\neq 0$ in $\Delta x=\sum
x^{\prime}\otimes x^{\prime\prime}.$

\section{Formal $\infty$-implication sequences}

Let $x$ be an element of a Hopf algebra over a finite field. In \cite{browder},
W. Browder introduced the notion of $\infty$-implications (of an infinite sequence)
associated with $x$ in the Hopf algebra. The following can be thought of as a
generalization of this: Let $x^{\smile_{\mathbf{1}}p}$ denote
the (right most) $p^{th}$-power of $x$ with respect to $\smile_{\mathbf{1}}
$-product with the convention that $x^{\smile_{\mathbf{1}}1}=x.$

\begin{definition}
\label{fis} Let $x\in V^{k},\,k\geq2,$ $d_{h}x\in\mathcal{D}_{\Bbbk}.$ A
sequence $\mathbf{x}=\{x(i)\}_{i\geq0}$  is a \underline{formal
$\infty$-implication sequence (f.i.s.)} of $x$ if

\begin{enumerate}
\item[\textit{(i)}] $x(0)=x,\, x(i)\in V^{(i+1)k-i},$ and
 $x(i)$ is not $\mu$-homologous to zero for all $i;$

\item[\textit{(ii)}] Either $x(i)=x^{\smile_{\mathbf{1}}(i+1)}$ or $x(i)$ is
resolved from the following relation in the filtered hga model $(RH,d_{h}):$
\begin{equation}\label{relation}
d_{h}{\mathfrak{b}}(i)=x^{\smile_{\mathbf{1}}(i+1)}+z(i)+\mu'{x}(i),\ \ \
{\mathfrak{b}}(i)\in V,\,z(i)\in\mathcal{D},\ \mu' \ \text{is divisible by}\ \mu.
\end{equation}

\end{enumerate}
\end{definition}

We are interested in the existence of an f.i.s. for an odd dimensional $x\in V.$

\begin{proposition}
\label{existence} Let $x\in V$ be of odd degree with $d_{h}x\in\mathcal{D} _{\Bbbk}$
such that ${x}$ is not $\mu$-homologous to zero. For $\mu\geq2,$ assume, in addition,
  there is no relation $d_{h}u=\mu x\!\mod {\mathcal D},$ some $u\in V.$
Then $x$ has an f.i.s. $\mathbf{x}=\{x(i)\}_{i\geq 0}.$
\end{proposition}

\begin{proof}
Suppose we have constructed $x(i)$ for $0\leq i< n.$ If $x ^{\smile_{\mathbf{1}}(n+1)}$
is not $\mu$-homologous to zero, set $x(n) =x^{\smile_{\mathbf{1}}(n+1)};$ otherwise,
there is the relation $d_{h}u=x^{\smile_{\mathbf{1}}(n+1)}+z+\mu' v$ for some $u,v\in
V,$ $z\in\mathcal{D}$ and $\mu'$ divisible by $\mu.$ Using (\ref{dif})--(\ref{prod})
one can easily establish the fact that $dx^{\smile_{\mathbf{1}}(n+1)}$ contains a
summand component of the form $-\sum_{k+\ell=n+1}\binom{n+1}{k}x^{\smile_{1}k}x^{\smile
_{1}\ell},\,k,\ell\geq1.$  We have that $v\neq 0$ in the aforementioned relation since
Proposition \ref{basic} (applied for $c=x^{\smile_{\mathbf{1}}(n+1)}$ and $a\cdot
b=-\binom{n+1}{k}x^{\smile_{1}k}\cdot x^{\smile _{1}\ell},$ some $k$). Clearly,
$d_{h}v=-\frac{1}{\mu'}d_h\left(x^{\smile_{\mathbf{1}}(n+1)}+z\right)\in\mathcal{D};$
Assuming $\mu'$ to be maximal  $v$
 is not $\lambda$-homologous to zero.
 Set $x(n)=v$
and $\mathfrak{b}(n)=u,\, z(n)=z$ to obtain (\ref{relation}) for $i=n.$
\end{proof}
\newpage
Thus, for $\mu=0$ (when $\Bbbk$ is a field of characteristic zero, for example)
 $\mathbf{x}=\{ x^{\smile_{\mathbf{1}}(n+1)} \}_{n\geq0}.$

\begin{remark}
1. The restriction on $x$ in Proposition \ref{existence} that no relation $d_{h}u=\mu
x\!\mod {\mathcal D}$ exists is essential. A counterexample is provided by the
exceptional group $F_{4}$: Let $A=C^{\ast}(BF_{4}; \mathbb{Z}_{3})$ be the cochain
complex of the classifying space $BF_{4}.$ Then we have the relation $du=3x$ in
$(RH,d)$  corresponding to the Bockstein cohomology homomorphism $\delta x_{8}=x_{9}$
on $H^{\ast}(BF_{4};\mathbb{Z}_{3})$ (in the notation of \cite{Toda}), but the element
$x(2)$ does not exist (see \cite{sane2} for more details).

2. Note that if $du=\mu x$ in Proposition \ref{existence}, but $[u][x]\neq 0\in
H_{\Bbbk},$ then   one can modify the proof of the proposition to show that $x$ again
has an f.i.s. $\{x(i)\}_{i\geq 0}.$  Note that in the above example we just have
$[u][x]=0\in H_{\mathbb{Z}_{3}}= H^{\ast}(BF_{4};\mathbb{Z}_{3}).$

3. The existence of $\infty$-implications of $x$ in \cite{browder} uses both the
$\smile$-product and the Pontrjagin product in the loop space (co)homology. In our case
each component of the sequence $\mathbf{x}$ is determined by item (ii) of
Definition \ref{fis} in which the first case can be thought of as related to the
$\smile$-product, and the second with the Pontrjagin product. In particular,
primitivity of $x$ required in \cite{browder} is not issue for  the existence of
$\infty$-implications of $x.$
\end{remark}

In certain cases, a given odd dimensional $b\in V$ rises to an infinite sequence
$\mathbf{b}=\{b_{i}\}_{i\geq0}$ with $b=b_{0}$ in the hga resolution $(RH,d).$ These
sequences are built by explicit formulas and include also the case $du=\lambda b,$
i.e., when the hypothesis of Proposition \ref{existence} formally fails (see, for
example, Case I of the proof of Proposition \ref{one} below). Namely, we have the
following cases:

(i) For $b\in V^{0,*}$ and $[b]^{2}=0\in H$ (i.e., there exists $b_{1}\in V^{-1,*}$
with $db_{1}=b^{2};$ e.g. $b_{1}=ab+\frac{\lambda-1}{2}b\smile_{1} b$ for $da=\lambda
b$ with $\lambda$ odd, some $a\in V^{-1,*}$), $\mathbf{b}=\{b_{i}\}_{i\geq0}$ is given
by
\begin{equation}
\label{xrelation1}db_{n}=\sum_{i+j=n-1} b_{i}b_{j}%
\end{equation}
and satisfies the following relation with ${\mathfrak{c}}_{i}\in V$
\[
d{\mathfrak{c}}_{n}=-(-1)^{n}\!\left( (n+1)b_{n}+b_{0}\!\smile_{1}\!
b_{n-1}\right) +\!\!\!\sum_{i+j=n-1}\!\!(-1)^{i}\left( {\mathfrak{c}}_{j}%
b_{i}- b_{i}{\mathfrak{c}}_{j} \right) ,n\geq1;
\]

(ii) For $b\in V^{0,*}$ and $[b]^{2}\neq0\in H$ (and $b_{1}=b\smile_{1} b$),
$\mathbf{b}=\{b_{i}\}_{i\geq0}$ is given by
\begin{equation}
\label{xrelation2}db_{2k}=\sum_{i+j=2k-1} b_{i}b_{j},\ \ \ \ \ \ \ db_{2k+1}%
=\sum_{i+j=k}( 2 b_{2i}b_{2j}+b_{2i-1}b_{2j+1}),
\end{equation}
and satisfies the following relation with ${\mathfrak{c}}_{i}\in V$ (below
${\mathfrak{c}}_{1}=0$)
\begin{multline*}
d{\mathfrak{c}}_{2k}=-(2k+1)b_{2k}-b_{0}\smile_{1} b_{2k-1}+ \sum _{i+j=k}2\left(
{\mathfrak{c}}_{2j-1}b_{2i}-b_{2i}{\mathfrak{c}}_{2j-1}\right)
\\
\hspace{2.25in}- \sum_{i+j=k}\left( {\mathfrak{c}}_{2j}b_{2i-1}-b_{2i-1}%
{\mathfrak{c}}_{2j}\right) ,\\
d{\mathfrak{c}}_{2k+1}=(k+1)b_{2k+1}+b_{0}\smile_{1} b_{2k}+\sum
_{i+j=2k}(-1)^{i}\left( {\mathfrak{c}}_{j}b_{i}-b_{i}{\mathfrak{c}}_{j} \right)
,\,k\geq1;
\end{multline*}

(iii) For $b\in V^{-1,*}$ and $db=\mu c,\, \mu\geq2,\, c\in V^{0,*} $ (below
$\omega_{0}:=c$), $\mathbf{b}=\{b_{i}\}_{i\geq0}$ is given by
\begin{multline}
\label{xrelation3}db_{n}=\sum_{i+j=n-1} b_{i}b_{j} +\mu c_{n},\\
c_{n}=-\omega_{0}\smile_{1}b_{n-1}- \sum_{\substack{i+j=n-1 \\i\geq1;\,j\geq 0}}
(-1)^{i}\omega_{i}\smile_{1} b_{j}-(-1)^{n}\omega_{n},\,n\geq1
\end{multline}
and satisfies the following relation with ${\mathfrak{c}}_{i}\in V$
\begin{multline*}
d{\mathfrak{c}}_{1}=2b_{1}+b_{0}\smile_{1} b_{0}+ \mu\omega_{0}\cup_{2}
b_{0},\\
d{\mathfrak{c}}_{n}=-(-1)^{n}\!\left( (n+1)b_{n}+b_{0}\smile_{1}
b_{n-1}\right)  +\sum_{i+j=n-1}\!\!(-1)^{i}\left( {\mathfrak{c}}_{j}%
b_{i}-b_{i}{\mathfrak{c}}_{j} \right) \\
\hspace{4.0in}+\mu\mathfrak{a}_{n},\\
\ \ \ \ \ \ \ \ \ \ \ \mathfrak{a}_{n}=\sum_{i+j=n-2}\!(-1)^{j}\left(
(\omega_{i}\!\cup_{2} b_{0})\smile_{1} b_{j} + \omega_{i}\!\smile_{1}
{\mathfrak{c}}_{j+1} \right) \! +\omega_{n-1}\!\cup_{2} b_{0},\\
\\
d\omega_{k}=\sum_{i+j=k-1}\mu\omega_{i}\smile_{1}\omega_{j},\ \ \ \omega
_{k}=\mu^{k}\omega_{0}^{\cup_{2}(k+1)},\,k\geq1,\,n\geq2.
\end{multline*}

\vspace{0.1in}

For example, in view of Proposition \ref{barVX}, the formulas above are enough to
calculate the loop space cohomology algebra with coefficients in $\Bbbk$ for Moore
spaces, i.e., the $CW$-complexes obtained by attaching an $(n+1)$-cell to the
$n$-sphere $S^{n}$ by a map $S^{n}\rightarrow S^{n}$ of degree $\mu.$

\subsection{Odd dimensional element $l(a)$}
\label{even}

 Given $m\geq2,$ let $H(A)$ be finitely generated as a $\Bbbk $-module with
$H^{i}(A)=0$ for $i>m.$  Let
 ${\mathcal{Z}}_{\Bbbk}$ be the subset of $RH$ defined by
\[\mathcal{Z}_{\Bbbk}=\mathcal{Z}'_{\Bbbk}+ \mathcal{Z}''_{\Bbbk}+\mathcal{D}_{\Bbbk},\]
\[
\mathcal{Z}'_{\Bbbk}=\left\{  v\in V\,|\,du=\lambda v,\
\ u\in V,\,\lambda\in\mathbb{Z}_{\Bbbk }\right\}  \] and
\[
\mathcal{Z}''_{\Bbbk}=\left\{  v\in V\,|\,   v=\lambda u,\ \ u\in {V},\, \lambda\in\mathbb{Z}\setminus\mathbb{Z}_{\Bbbk}  \right\} .
\]
Given $x\in V$ with $d_{h}x=w\in{\mathcal{Z}}_{\Bbbk},$ $ w=w'+w''+z,$  define
\[
\tilde{x}=
 \frac{l.c.m.(\lambda'';\mu)}{\lambda''}(\lambda' x-u), \, du=\lambda'w',\,
w''=\lambda'' v'',
\]
to obtain $d_{h}\tilde{x}\in{\mathcal{D}}_{\Bbbk}.$

Regarding (\ref{general}), define also the following subsets
$K^*_{\mu},K^*_{0}\subset\mathcal{V}^{-1,\ast}$ with $K^*_{\mu}\subset K^*$ as
\[
\begin{array}
[c]{lll} \ K_{\mu}=\left\{  a\in {K}\,|\,
\lambda\ \text{is divisible by}\ \mu  \right\}  ,\vspace{1mm} &  & \\
K_{0}=\left\{  u\in\mathcal{V}^{-1,\ast}\setminus\mathcal{E}\,|\,du\in
\mathcal{D}^{0,\ast}\right\},

\end{array}
\]
and assign to a given even dimensional element $a\in V^{0,\ast}\cup K_{\mu} $
an odd dimensional element $l(a)\in V$ with $dl(a)\in {\mathcal D}_{\Bbbk}$ as
follows. If $a\in V^{0,\ast},$ let $l(a)\in K_0$ be an element such that
$dl(a)=a^{k},$ where $k\geq2$ is chosen to be the smallest. If $a\in K_{\mu}$ with
$
da=\lambda b
$
consider the relation
\begin{equation}\label{la1}
du_{1}=-a^{2}+\lambda v_{1},\ \ dv_{1}=\frac{1}{\lambda}d(a^{2}),\ \ u_{1}\in
V^{-3,\ast},\,v_{1}\in V^{-2,\ast},
\end{equation}
and the perturbation $hu_{1}=h^{2}u_{1}+h^{3}u_{1}.$ When $hu_{1}\in\mathcal{Z}_{\Bbbk},$   set
$l(a)=\widetilde{u}_{1},$ while when
   $h^3u_{1}\notin\mathcal{Z}_{\Bbbk},$ consider $\mathfrak{u}_1=h^3u_{1}|_{V^{0,*}},$ the component of $h^3u_{1}$ in $V^{0,*},$ and
 define
$l(a)$ as $l(\mathfrak{u}_1).$ When $h^2u_{1}\notin\mathcal{Z}_{\Bbbk},$ and  $h^3u_{1}\in\mathcal{Z}_{\Bbbk},$
 choose the smallest $n>1$
such that there is the relation
\begin{multline}\label{la2}
du_{n}=-ah^{2}u_{n-1}+\lambda v_{n},\ \ dv_{n}=\frac{1}{\lambda}d(a h^{2}
{u}_{n-1}),\ \ u_{n}\in V^{-3,\ast},\,v_{n}\in V^{-2,\ast},\\
\text{with}\ \ h^{2}u_{n}\in{\mathcal{Z}_{\Bbbk}}.
\end{multline}
(The inequality  $(n+1)|a|>m$ guarantees the existence of  such a relation, since
$h^2u_i\in {\mathcal D}+K_{\mu},$ while $K^j_{\mu}=0$
for $j>m$  in the minimal $V\subset RH.$)
Then set
$l(a)=\widetilde{u}_{n}$ for $h^{3}u_{n}\in{\mathcal{Z}}_{\Bbbk};$ otherwise, define
$l(a)$ as $l(\mathfrak{u}_{n})$ for $\mathfrak{u}_n=h^3u_{n}|_{V^{0,*}}.$

\section{Proof of Theorem \ref{tau}}

\label{theorem} The proof of the theorem relies on the two basic propositions below in
which the condition that $\tilde{H}(A)$ has at least two algebra generators is treated
in two specific cases.

\begin{proposition}
\label{one} Let $H_{\Bbbk}$ be a finitely generated $\Bbbk$-module with $\mu\geq 2.$ If
$\tilde{H}_{\Bbbk}$ has at least two algebra generators and $\tilde {H}_{\mathbb{Q}}$
is either trivial or has a single algebra generator, there are two sequences of odd
degree elements $\mathbf{x}_{\Bbbk}=\{x(i) \}_{i\geq0}$ and
$\mathbf{y}_{\Bbbk}=\{y(j)\}_{j\geq0}$ in $V_{\Bbbk}$ whose degrees form arithmetic
progressions such that all $\bar{x}(i),$ $\bar{y}(j)$ are $\bar{d}_{h}$-cocycles in
$\bar{V}_{\Bbbk}$ and the classes $\left\{
[s^{-1}(x(i)\smile_{\mathbf{1}}y(j))]\right\}  _{i,j\geq0}$ are linearly independent in
$H(\bar{V}_{\Bbbk},\bar{d}_{h}).$
\end{proposition}

\begin{proof}
The hypotheses of the proposition imply that $K_{\mu}$  defined in subsection \ref{even}  above is non-empty; also by the
restriction on $\tilde{H}_{\mathbb{Q}},$ relation (\ref{general}) reduces to
\[
da=\lambda b^{m},\,\lambda\neq0,\,m\geq1,\,(\lambda,m)\neq(1,1),\,b\in
\mathcal{V}^{0,\ast}%
\]
for $a\in{\mathcal{V}}^{-1,\ast}$ to be of the smallest degree.

In the three cases below, we exhibit two odd dimensional elements $x,y\in V\setminus
{\mathcal E}$ that fail to be $\mu$-homologous to zero. \vspace{1mm}

Case I. Let $a\in K_{\mu}$ be of the smallest degree in $K_{\mu}\cup K_{0}$ with $da=\lambda
b$ and let $|a|$ be even.
 Consider the element $l(a).$ If it is not $\lambda$-homologous to zero, set $x=l(a);$
  otherwise,
 we must have  relation
(\ref{general}) in which $v_i=a$ for some $i$ and
 $hu\in \mathcal{Z}_{\Bbbk}$ with $|u|<|l(a)|,\, u\in \bigcup_{i\geq 1}\mathcal{V}^{-i,*}\setminus {\mathcal E}. $
 By (\ref{general}) choose     $u$ to be of the smallest degree  with $hu\in \mathcal{Z}_{\Bbbk},\, u\neq
u_i,a_1,$ where $u_i$ is given by (\ref{la1})--(\ref{la2}) and $da_1=-ab+\lambda b_1,\, db_1=b^2.$
  Set $x=\tilde{u}$ for $|u|$  odd. If $|u|$ is even and $u\in \bigcup_{i> 1}\mathcal{V}^{-i,*}\setminus {\mathcal E}$
set $x=\tilde{v};$ if $u\in K_0$
 and $du$ contains an odd dimensional $v_i\in V^{0,*}$ with $[v_i]\neq 0\in H_{\mathbb Q},$
set $x=v_i;$ otherwise,
 for each
monomial $P_{s}(v_{1},...,v_{r_{s}})$ choose a variable $v_{i}$ with a
relation $du_{i}=\mu_{i}v_{i}$ (for example, we can choose  $v_i$ to be odd dimensional for all $s$). Let $\lambda$ be the smallest integer
divisible by all $\mu_{i},$ and replace $v_{i}$ by $\frac{\lambda}{\mu_{i}
}u_{i}$ to detect a new relation in $(RH,d)$  given again by (\ref{general}):
\[
d{w}=\sum_{1\leq s\leq n}\frac{\lambda_{s}\lambda}{\mu_{i}}P_{s}(v_{1},...,v_{i-1},
u_{i},v_{i+1}...,v_{r_{s}})+\lambda u,\ \ \ \ \lambda_s\in \mathbb{Z}_{\Bbbk},\  w\in
\mathcal{V}^{-2,\ast}.
\]
Hence, $|w|$ is odd, and set $x=\widetilde{w}$ for $h^{2}w
\in{\mathcal{Z}}_{\Bbbk}.$ If $h^{2}w\notin{\mathcal{Z}}_{\Bbbk}$ we have
the following two subcases:

(i1) Assume there exists $v\in K_{\mu}$ with $dv=\lambda h^{2}w.$
If
 ${[\bar{v}]}_{\lambda}\neq0,$  set $x=v;$ otherwise we have a relation
 $d_{h}u'=v+z+\lambda^{\prime}v^{\prime},$ some $\,u',v^{\prime}\in
V,\,z\in\mathcal{D}.$ Clearly, $h^{tr}v'=-\frac{\lambda}{\lambda^{\prime}}h^2w\!
\mod{\mathcal D} ,$ and set $x=\frac{\lambda}{\lambda^{\prime}}w+v^{\prime}.$ Note that
$x$ is not $\lambda$-homologous to zero since the component
$\frac{\lambda^2}{\lambda^{\prime}}u$ in $dx.$

(i2) Assume $[h^{2}w]\neq0\in H_{\mathbb{Q}}.$ When $r_s>1 $ for all s,
  choose  a variable $v_j$ different from  $v_i$ in $P_{s}(v_{1},...,v_{r_{s}})$
to form $w'$ entirely analogously to $w,$ and then  find $x$ similarly to the above
unless $[h^{2}w']\neq 0\in H_{\mathbb{Q}},$ in which case set $x=\alpha w+\beta w',$
some $\alpha,\beta \in {\mathbb Z}.$ When $\mathbf{k}= \{s\in \underline{n} \,|\, r_s=1
\ \text{in}\ du\}\neq \varnothing ,$  i.e.,
 $P_{s}(v_1,...,v_{r_{s}})= v_1^{2m_1+1}:= v_s^{2m_s+1},  m_s\geq 1,$
 $|v_s|$ is odd  for $s\in \mathbf{k}$ (in particular, $\mu_s $ is even, since
$[v_s]^2=0\in H$
 for  $\mu_s$ odd; c.f. (\ref{xrelation1})), then
\[ du'=\left\{ \begin{array}{llll}
\sum_{s\in \mathbf{k}} \frac{\lambda_{s}\lambda}{2}(v_{s}\smile_1
v_{s})v_{s}^{2m_{s}-1}\vspace{1mm}\\+ \sum_{s\in \underline{n}\setminus
\mathbf{k}}\frac{\lambda_{s}\lambda}{\mu_{j}}P_{s}(v_{1},...,v_{j-1},
u_{j},v_{j+1}...,v_{r_{s}})+\lambda u , &  \mathbf{k}\neq \underline{n},\vspace{0.1in}\\

\sum_{s\in \underline{n}} \lambda_s(v_{s}\smile_1v_{s})v_{s}^{2m_{s}-1}+2u,&
\mathbf{k}=\underline{n}

 \end{array}
 \right.
 \]
 with
 $u'\in {\mathcal V}^{-2,\ast},$
  and by considering $h^2u'$ we find $x$ as in item (i1).

To find $y,$ consider $b$ and the associated sequence $\mathbf{b}=\{b_{i}\}$ given by
(\ref{xrelation1}) or (\ref{xrelation2}). If $hb_{i}\in{\mathcal{Z} }_{\Bbbk}$ for all
$i,$ set $y=b$ and $\mathbf{y}=\{\tilde{b}_{i}\}_{i\geq0}.$ If
$h\mathbf{b}\nsubseteq{\mathcal{Z}}_{\Bbbk},$ consider the smallest $p>0$ such that
$h^{tr}b_{p}\notin{\mathcal{Z}}_{\Bbbk}.$  Consider $t_p=h^{tr}b_{p}|_{V^{0,*}},$  and
if $\left[\overline{l(t_p)}\right]_{\lambda}\neq 0,$ set $y=l(t_p);$ if
$\left[\overline{l(t_p)}\right]_{\lambda}=0$  and
 $\alpha h^{3}u_{i}+\beta
h^{tr}b_{p}=0,\, \alpha,\beta\in\mathbb{Z},$ for some $u_{i}$ from (\ref{la1})--(\ref{la2}),
set $y=\alpha u_{i}+\beta b_{p};$ otherwise, we obtain
$l(t_p)\in K_0$ different from $l(a)$ above; consequently,
 we must have  another relation in $(RH,d)$ given by (\ref{general}) in which $v_i=t_{p} $ for
some $i$ and
 $hu\in \mathcal{Z}_{\Bbbk}$ with $|u|<|l(t_p)|,$
and then
 $y$ is found  similarly to $x.$

\vspace{1mm} Case II. Let $a\in K_{\mu}$ be of the smallest degree in $K_{\mu}\cup K_{0}$
with $da=\lambda b$ and let $|a|$ be odd. Set $x=a.$ Consider $l(b)\in K_{0},$ and then
$y$ is found as in Case I.

\vspace{1mm}

Case III. Let ${a}\in K_{0}$ be of smallest degree in $K_{\mu}\cup K_{0}$ with
$da=\lambda b^{m},m\geq 2,$    and $[b]\neq 0\in
H_{\mathbb{Q}}$.  Set
\[
x=\left\{
\begin{array}
[c]{llll}%
b, & |b| & \text{is odd} & \\
a, & |b| & \text{is even}. &
\end{array}
\right.
\]
To find $y$ consider the following two subcases:

(i) Assume $\lambda\in\mathbb{Z}\setminus\mathbb{Z}_{\Bbbk}.$ When both $|a|$ and $|b|$
are odd, set $y=a;$ otherwise, either $|a|$ or $|b|$ is even,      in which case
consider $l(\tilde{a})$ or $l(b)$ respectively, and then
 $y$ is
found  as in Case I.

(ii) Assume $\lambda\in\mathbb{Z}_{\Bbbk}.$ Since $K_{\mu}\neq\varnothing,$  this subcase reduces either to  Case I or to Case II.

Finally, having found the elements $x$ and $y$ in Cases I-III, consider the f.i.s.
$\mathbf{x}$ and $\mathbf{y}$ in $V$ and the induced sequences
$\mathbf{x}_{\Bbbk}=\{x(i)\}_{i\geq 0}$ and $\mathbf{y}_{\Bbbk}=\{y(j)\}_{j\geq 0}$ in
$V_{\Bbbk}.$ Then the both sequences $\bar{\mathbf{x}}_{\Bbbk}$ and
$\bar{\mathbf{y}}_{\Bbbk}$ consist of $\bar{d}_{h}$-cocycles in $\bar{V}_{\Bbbk}$ whose
degrees form an arithmetic progression respectively. Thus, we obtain that
$[\bar{\mathbf{x}}_{\Bbbk}],\,[\bar{\mathbf{y}}_{\Bbbk}]\subset H(\bar
{V}_{\Bbbk},\bar{d}_{h})$ are sequences of non-trivial classes. Moreover, they are
linearly independent and $\left\{  [s^{-1}(x(i)\smile_{\mathbf{1}} y(j))]\right\}
_{i,j\geq0}$ is the sequence of linearly independent classes in
$H(\bar{V}_{\Bbbk},\bar{d}_{h})$ as required.

\vspace{1mm}
\end{proof}

Before proving the second basic proposition we need the following auxiliary statement.
Given a cochain complex $(C^{\ast},d)$ over $\mathbb{Q},$ let
$S_{C}(T)=\sum_{n\geq0}(\dim_{\mathbb{Q}}C^{n})T^{n}$ and $S_{H(C)}
(T)=\sum_{n\geq0}(\dim_{\mathbb{Q}}H^{n}(C))T^{n}$ be the Poincar\'{e} series. As
usual, we write $\sum_{n\geq0}a_{n}T^{n}\leq\sum_{n\geq0}b_{n}T^{n}$ if and only if
$a_{n}\leq b_{n}.$ The following proposition can be thought of as a modification of
Propositions 3 and 4 in \cite{Sul-Vig} for the non-commutative case.

\begin{proposition}
\label{quotient} Given an element $y\in V_{_{\mathbb{Q}}}$ of total degree
$K_{\mu}\geq2$ such that $\bar{d}_{h}(\bar{y})=0,$ let $y\bar{V}_{_{\mathbb{Q}}
}\subset\bar{V}_{_{\mathbb{Q}}}$ be a subcomplex (additively) generated by the
expressions $\{\bar{y}=s^{-1}y,\,s^{-1}(y\smile_{1} v)\}_{v\in V_{_{\mathbb{Q} }}}.$
Then
\begin{equation}
\label{inequality}S_{H(\bar{V}_{_{\mathbb{Q}}}/y\bar{V}_{_{\mathbb{Q}}})}
(T)\leq(1+T^{k-1})S_{H(\bar{V}_{_{\mathbb{Q}}})}(T).
\end{equation}

\end{proposition}

\begin{proof}
Consider the inclusion of cochain complexes $s^{k}\bar{V}_{_{\mathbb{Q}}
}\overset{\iota}{\rightarrow}\bar{V}_{_{\mathbb{Q}}}$ defined for
$1\in\mathbb{Q}=(s^{k}\bar{V}_{_{\mathbb{Q}}})^{k}$ by $\iota(1)=\bar{y},$ and for
$s^{k}(\bar{v})\in(s^{k}\bar{V}_{_{\mathbb{Q}}})^{>k},\,v\in V_{_{\mathbb{Q}}}^{>1},$
by $\iota(s^{k}(\bar{v}))=s^{-1}(y\smile_{1}v).$ Then
$\iota(s^{k}\bar{V}_{_{\mathbb{Q}}})=y\bar{V}_{_{\mathbb{Q}}}$ and there is the short
exact sequence of cochain complexes
\[
0\rightarrow s^{k}\bar{V}_{_{\mathbb{Q}}}\overset{\iota}{\rightarrow}\bar
{V}_{_{\mathbb{Q}}}\rightarrow\bar{V}_{_{\mathbb{Q}}}/y\bar{V}_{_{\mathbb{Q}}
}\rightarrow0.
\]
Consider the induced long exact sequence
\[
\cdots\rightarrow H^{n-k}(\bar{V}_{_{\mathbb{Q}}})\overset{H^{n}
(\iota)}{\longrightarrow}H^{n}(\bar{V}_{_{\mathbb{Q}}})\rightarrow H^{n}
(\bar{V}_{_{\mathbb{Q}}}/y\bar{V}_{_{\mathbb{Q}}})\rightarrow H^{n-k+1}
(\bar{V}_{_{\mathbb{Q}}})\rightarrow\cdots.
\]
Let $I=\oplus I_{n},$ where $I_{n}=\operatorname{Im}(H^{n}(\iota)),$ $n\geq0,$ and form
the exact sequence
\[
0\rightarrow I_{n}\rightarrow H^{n}(\bar{V}_{_{\mathbb{Q}}})\rightarrow
H^{n}(\bar{V}_{_{\mathbb{Q}}}/y\bar{V}_{_{\mathbb{Q}}})\rightarrow
H^{n-k+1}(\bar{V}_{_{\mathbb{Q}}})\rightarrow I_{n+1}\rightarrow0.
\]
Since $I_{0}=0,$ we have
\[
\sum_{n\geq0}(\dim_{\mathbb{Q}}I_{n}+\dim_{\mathbb{Q}}I_{n+1})T^{n}
=\frac{(1+T)S_{I}(T)}{T}.
\]
Now apply the Euler-Poincar\'{e} lemma for the above exact sequence to obtain the
equality
\[
\frac{(1+T)S_{I}(T)}{T}-S_{H(\bar{V}_{_{\mathbb{Q}}})}(T)+S_{H(\bar
{V}_{_{\mathbb{Q}}}/y\bar{V}_{_{\mathbb{Q}}})}(T)-T^{k-1}S_{H(\Bar
{V}_{_{\mathbb{Q}}})}(T)=0.
\]
Consequently,
\[
S_{H(\bar{V}_{_{\mathbb{Q}}}/y\bar{V}_{_{\mathbb{Q}}})}(T)=(1+T^{k-1}
)S_{H(\bar{V}_{_{\mathbb{Q}}})}(T)-\frac{(1+T)S_{I}(T)}{T},
\]
and since $S_{I}(T)\geq0,$ we get (\ref{inequality}) as required.
\end{proof}

\begin{proposition}
\label{rationalone} Let $H_{\Bbbk}$ be a finitely generated $\Bbbk$-module. If
$\tilde{H}_{_{\mathbb{Q}}}$ has at least two algebra generators and
$A_{_{\mathbb{Q}}}=A^{\prime}\otimes_{\mathbb{Z}}\mathbb{Q},$ the set $\left\{
\tau_{i}(A_{_{\mathbb{Q}}})=\dim_{\mathbb{Q}}Tor_{i}
^{A_{_{\mathbb{Q}}}}(\mathbb{Q}\,,\mathbb{Q})\right\}  $ is unbounded.
\end{proposition}

\begin{proof}
Consider the first two generators $a_{i}\in V^{-1,*}_{_{\mathbb{Q}}}$ with
$da_{i}\in{\mathcal{D}^{0,*}},$ $i=1,2.$ We have two cases:

(i) Both $|a_{1}|$ and $|a_{2}|$ are odd. Set $x=a_{1}$ and $y=a_{2}.$ Then both
$\bar{x}$ and $\bar{y}$ are $\bar{d}_{h}$-cocycles and the classes $[\bar{x}]$ and
$[\bar{y}]$ are non-trivial in $H(\bar{V}_{_{\mathbb{Q}}} ,\bar{d}_{h}).$ Consequently,
the classes
\begin{equation}
\left\{  \lbrack s^{-1}\left(  x^{\smallsmile_{\mathbf{1}}i}\smile
_{\mathbf{1}}y^{\smallsmile_{\mathbf{1}}j}\right)  ]\right\}  _{i,j\geq
1}\label{classes}%
\end{equation}
are linearly independent in $H(\bar{V}_{_{\mathbb{Q}}},\bar{d}_{h}).$

(ii) Either $|a_{1}|$ or $|a_{2}|$ is even. Denote the (smallest) even dimensional
generator by $a$ and consider $da.$ Then for $a,$ (\ref{general}) reduces to
\[
da=uv,\ \ u\in V_{_{\mathbb{Q}}}^{0,2k+1}\ \ \text{and} \ \ v\in
R^{0}H_{_{\mathbb{Q}}}^{2\ell },\ \text{some}\ k,\ell\geq1.
\]
There are the following induced relations in $(RH_{_{\mathbb{Q}}},d):$
\[%
\begin{array}
[c]{llll} db=-u(a+u\smile_{1}v)-au, & b\in V_{_{\mathbb{Q}}}^{-2,2(2k+\ell
+1)}\ \ \ \ \text{and}\vspace{1mm}  \\
dc=-u\left(  v\smile_{1}a+(u\cup_{2}v)v+u(v\cup_{2}v)\right)  -a^{2}+bv, & c\in
V_{_{\mathbb{Q}}}^{-3,4(k+\ell)+2}.
\end{array}
\]
Thus we have $hc=h^{2}c+h^{3}c,$ and in particular, $dh^{2}c=h^{2}b\cdot v.$ Consider
the following two cases:

(1) Assume $hc\in{\mathcal{D}}.$ Set $x=u,$ $y=c,$ and obtain linearly independent
classes in $H(\bar{V}_{_{\mathbb{Q}}},\bar{d}_{h})$ by formula (\ref{classes}).

(2) Assume $hc\notin{\mathcal{D}}.$ Let $(\bar W,\bar {d}_{W})=(\bar
{V}_{_{\mathbb{Q}}}/\bar C,\bar{d}_{W}),$ where $C\subset {V}_{_{\mathbb{Q}}}$ is a
subcomplex (additively) generated by the expressions ${hc}$ and $ {hc\smile_{1}z}$ for
$z\in V_{_{\mathbb{Q}}}.$ Define $x$ and $y$ as the projections of the elements ${u}$
and ${c}$ from ${V} _{_{\mathbb{Q}}}$ under the quotient map
${V}_{_{\mathbb{Q}}}\rightarrow {V}_{_{\mathbb{Q}}}/C,$ respectively. Then $\bar x$ and
$\bar y$ are $\bar{d}_{W}$-cocycles in $\bar W.$ Once again apply formula
(\ref{classes}) to obtain linearly independent classes in $H(\bar W,\bar{d}_{W}).$
Finally, Proposition \ref{quotient} implies that $S_{H(\bar W)} (T)\leq
S_{H(V_{_{\mathbb{Q}}})}(T),$ and an application of Proposition \ref{barV} completes
the proof.
\end{proof}

\subsection{Proof of Theorem \ref{tau}}

In view of Proposition \ref{barV}, the proof reduces to the examination of the
$\Bbbk$-module $H(\bar{V}_{\Bbbk},\bar{d}_{h}).$ If $\tilde{H}_{\Bbbk}$ has a single
algebra generator $a,$ then the set $\left\{  \tau_{i}(A)\right\}  $ is bounded since
$\tau_{i}(A)=1.$ For example, this can be seen from the fact that
$H(\bar{V}_{\Bbbk},\bar{d}_{h})$ is generated by a single sequence induced by
(\ref{xrelation1}) or by (\ref{xrelation2}), where $x=a$ or $x=l(a)$ for $|a|$ odd or
even respectively, and by $\smile_{1}$-products of its components. If
$\tilde{H}_{\Bbbk}$ has at least two algebra generators, then the proof follows from
Propositions \ref{one} and \ref{rationalone}.

\vspace{5mm}

\end{document}